\documentclass[12pt]{amsart}
\usepackage{amssymb,latexsym,amsfonts, amsmath}
\usepackage{hyperref}
\usepackage{lineno}
\usepackage{enumitem}
\usepackage{chngcntr}
\usepackage{pst-all}
\begin{document}
\newtheorem{theorem}{Theorem}[section]
\newtheorem{definition}[theorem]{Definition}
\newtheorem{proposition}[theorem]{Proposition}
\newtheorem{lemma}[theorem]{Lemma}
\newtheorem{remark}[theorem]{Remark}
\newtheorem{corollary}[theorem]{Corollary}
\newtheorem{question}{Question}
\newtheorem{example}{Examples}[section]
\newtheorem{notation}[theorem]{Notation}
\newtheorem{claim}[theorem]{Claim}
\newtheorem{fact}[theorem]{Fact}
\newcommand\cl{\begin{claim}}
\newcommand\ecl{\end{claim}}
\newcommand\rem{\begin{remark}\upshape}
\newcommand\erem{\end{remark}}
\newcommand\ex{\begin{example}\upshape}
\newcommand\eex{\end{example}}
\newcommand\nota{\begin{notation}\upshape}
\newcommand\enota{\end{notation}}
\newcommand\dfn{\begin{definition}\upshape}
\newcommand\edfn{\end{definition}}
\newcommand\cor{\begin{corollary}}
\newcommand\ecor{\end{corollary}}
\newcommand\thm{\begin{theorem}}
\newcommand\ethm{\end{theorem}}
\newcommand\prop{\begin{proposition}}
\newcommand\eprop{\end{proposition}}
\newcommand\lem{\begin{lemma}}
\newcommand\elem{\end{lemma}}
\newcommand\fct{\begin{fact}}
\newcommand\efct{\end{fact}}
\providecommand\qed{\hfill$\quad\Box$}
\newcommand\pr{{Proof:\;}}
\newcommand\prc{\par\noindent{\em Proof of Claim: }}
\def\R{{ \mathbb R}}
\def\Z{{\mathbb Z}}
\def\N{{\mathbb N}}
\def\F{{\mathbb F}}
\def\C{{\mathbb C}}
\def\si{{\sigma}}
\def\Ps{{\mathcal P}}
\def\Ra{{\mathcal R}}
\def\L{{\mathcal L}}
\def\M{{\mathcal M}}
\def\Ca{{\mathcal C}}
\def\U{{\mathcal U}}
\def\I{{\mathcal I}}
\def\W{{\mathcal W}}
\def\S{{\mathcal S}}
\def\B{{\mathfrak B}}
\def\FF{{\mathfrak F}}
\def\V{{\mathcal V}}
\def\K{{\mathcal K}}
\def\O{{\mathcal O}}
\def\T{{\mathcal T}}
\def\A{{\mathcal A}}
\def\cF{{\mathcal F}}
\def\D{{\mathcal D}}
\def\ki{{\bar K}}
\def\G{{\mathcal G}}
\def\l{{\ell}}
\def\bc{{\bold{c}}}
\author{Fran\c coise {Point}}
\address{Department of Mathematics (De Vinci)\\ UMons\\ 20, place du Parc 7000 Mons, Belgium}
\email{point@math.univ-paris-diderot.fr}
\title{A decidable expansion of $(\Gamma,+,F)$ with the independence property.}
\begin{abstract}  
Let $(\Gamma,+,F)$ be a finitely generated $\Z[F]$-module where $F$ is an injective endomorphism of the abelian group $\Gamma$. We restrict ourselves to a finite automa presentable subclass, introduced by J. Bell and R. Moosa \cite{BM} and define an expansion containing the $\cF$-sets defined by R. Moosa and T. Scanlon \cite{MS},  where every automatic subset is definable.
\end{abstract}
\date{\today}
\maketitle
\section{Introduction}
\par Let $(\Gamma,+, F)$ be a finitely generated $\Z[X]$-module where $F$ be an endomorphism of $(\Gamma,+,0)$ and the action of $X$ is given by the action of $F$; we will denote these modules $\Z[F]$-modules. We will always assume that $F$ is injective.
\par Under the further assumption that these modules are finitely generated abelian groups, the class of finitely generated $\Z[F]$-modules we consider, has been introduced by R. Moosa and T. Scanlon in their work on the Mordell-Lang conjecture for semi-abelian varieties over finite fields. They  enriched these finitely generated abelian groups $(\Gamma,+,F)$ by a collection of, so called, $F$-sets (for instance orbits by powers of $F$) (see Definition \ref{Fsets}), retaining the stability of their first-order theory. 
Their motivation was to understand the nature of subsets of the form $X\cap \Gamma$, where $\Gamma$ is a finitely generated subgroup of a semiabelian variety $V$ over a finite field $\F_{q}$, $X$ is a closed subvariety of $V$ and $\Gamma$ is invariant under the action of $q$-power Frobenius endomorphism $F$ \cite{MS}.
\par Later J. Bell and R. Moosa showed that these $F$-sets $\cF$ are automatic \cite[Theorem 6.9]{BM} and this aspect was further investigated by C. Hawthorne \cite{H}. 
 \par We aim to define an expansion of $(\Gamma,+,F,\cF)$ which is finite-automaton (FA) presentable where definable subsets coincide with automatic ones. This expansion will be automatically decidable but it will have 
 the independence property.
 \par  When $\Gamma$ is the group of integers and $F$ the multiplication by a prime number, one recovers a classical result of B\"uchi (see section \ref{buchi}). 
\par These finitely generated $\Z[F]$-modules have a finite generating set $\Sigma$  which is a $F^n$-spanning set, for some $n>0$, which allows them to be FA-presentable structures (see Definition \ref{span} and Fact \ref{rec}). When $\Gamma$ is a finitely generated abelian group, C. Hawthorne has shown that this is equivalent to require that the eigenvalues of $F\otimes (1\restriction \C)$ acting on the $\C$-subspace $\Gamma\otimes \C$ are strictly bigger than $1$ (see section \ref{main}).
\par We will first review the notion of FA-presentable structures in section \ref{FA}, then placing ourselves in the special case where $\Gamma=\Z$, we will recall the result of B\"uchi (section \ref{buchi}). Then we will describe some of the previous results on this subclass of finitely generated $\Z[F]$-modules  $(\Gamma,+,F)$ (section \ref{main}).
Finally we will define a sufficiently rich expansion of $(\Gamma,+,F)$ in order to use the fact that automatic sets coincide with regular languages \cite{PP}. 
 \section{Finite-automaton presentable structures}\label{FA}
 In order to be self-contained we briefly review below the notion of {\it automatic} structures. We begin by recalling the definition of finite automata. 

 \subsection{Finite automata}\label{aut}
 There are several equivalent definitions of finite automata. A classical reference is the book of S. Eilenberg \cite{E} (see also \cite[chapter 1]{PP}).
   \par  A finite automaton $\A$ is a finite state machine given by the following data $(\Sigma, Q, q_0, F_a, T)$, where $\Sigma$ is a finite alphabet, $q_0$ the initial state, $Q$ a finite set of states, $F_a\subset Q$ a set of accepting states and $T$ a transition function from $Q\times \Sigma\to Q$ with the convention that $T(q,\lambda)=q$ for $\lambda$ the empty word.
  \par Let  $\Sigma^*$ be the set of all finite words on $\Sigma$. Let $\sigma\in \Sigma^*$ and $a\in \Sigma$, one extends the transition function to a function from $Q\times \Sigma^*\to Q$, still denoting it by $T$, by setting $T(q,\sigma a):=T(T(q,\sigma),a)$. 

 \par The automaton $\A$ accepts $\sigma\in \Sigma^*$ if $T(q_0,\sigma)\in F_a$. (We will say that $\A$ works on $\Sigma$.)

\par  A subset $L\subset \Sigma^*$ (or alternatively a language on $\Sigma$) is recognized/accepted by $\A$ if $L=\{\sigma\in\Sigma^*\colon T(q_0,\sigma)\in F_a\}$. 
\par Now let us recall the notion of regular languages on $\Sigma$. Given two words $\sigma, \tau \in \Sigma^*$, we denote their concatenation by $\sigma^{\frown}\tau$.
 \dfn
 The class $\mathcal R$ of {\it regular} languages is the smallest class of languages on $\Sigma$ that contains all finite languages and is closed under the following operations (called regular operations):
 if $L_0, L_1, L_2\in \mathcal R$, then 
 \begin{enumerate}
 \item $L_1\cup L_2\in \mathcal R$,
 \item $L_1^{\frown}L_2:=\{\sigma^{\frown}\tau\colon \sigma\in L_1, \tau\in L_2\}\in \mathcal R$,
 \item $L_0^*:=\{\sigma_1^{\frown}\ldots^{\frown}\sigma_n\colon n\in \N, \sigma_i\in L_0, 1\leq i\leq n\}\in \mathcal R.$
 \end{enumerate}
 \par A regular language $L$ is of {\it complexity $\leq c$} if $L$ can be obtained starting from $\{a\}$, $a\in \Sigma$ and $\{\lambda\}$, using $c$ times  these regular operations.
 \edfn
\par In the following, we are going to use Kleene's theorem on the correspondence between subsets of $\Sigma^*$ accepted (or recognized) by a finite automaton and regular languages on $\Sigma$. 
 \par {\bf Fact:} \cite[Chapter 1, Theorem 4.4]{PP} A subset $L\subset \Sigma^*$ is regular iff it is accepted by a finite automaton $\A$ working on $\Sigma$.
\subsection{Automatic structures}\label{automatic}
 Before recalling the notion of automatic structures, we need to recall the operation of convolution \cite[section 3]{Ho}. Given two elements $u, v\in \Sigma^*$ and an additional symbol $\sharp$, with $u=u_{1}\ldots u_{k}$ and $v=v_{1}\ldots v_{m}$, $k\leq m$,
 we define $u\star v$ as the word of length $m$ on $(\Sigma\cup\{\sharp\})^2$ of the form $(u_{1},v_{1})\ldots (u_{k},v_{k})(\sharp,v_{k+1})\ldots (\sharp,v_{m})$.
 \dfn \cite{Ho} Let $\L$ be a finite relational language. A first-order (countable) $\L$-structure $\M$ is automatic if there is a finite alphabet $\Sigma$ and an effective procedure to represent the domain of $\M$ by a regular language $D\subset \Sigma^*$ and using convolution, to represent the graph of each $n_{i}$-relation, $1\leq n_i$, $i\in I$, $I$ a finite set, by a regular 
 subset of $((\Sigma\cup\{\sharp\})^{n_{i}})^*$.
 \par When $\L$ is a finite but not a relational language, we will say that the countable $\L$-structure is automatic if the corresponding relational structure is automatic, namely we replace every function by its graph. (Of course it changes the notion of substructures but will not matter when considering decidability issues.)
 \edfn
\par It follows that in an automatic structure, whether an atomic formula holds in $\M$ can be checked by finite automata \cite[Lemma 3.6]{Ho} (which is the definition of an FA-presentable structure given in \cite[Section 2.1]{N}, where function symbols are allowed in $\L$). Note that being FA-presentable for a group or a ring is a strong restriction on its algebraic structure \cite[sections 3, 4]{N} (however any finitely generated abelian group is FA-presentable). This is reflected by the following result, which connects automaticity and decidability.
 (Hodgson attributes the paternity of that result to R. B\"uchi and C. Elgot).
 \vskip 0,3cm
 \fct \cite[Theorem 3.5]{Ho} If a countable structure $\M$ is automatic, then it is decidable.
 \efct
The proof of the above result uses the closure of the class of finite automata by the operations corresponding to the boolean operations on formulas and the existential quantification. One then deduces that any definable subset in $\M$ is recognized by a finite automaton. The result follows from decidability of  
 the emptiness problem for finite automata.

 \section{Definability and finite automata in the natural numbers}\label{buchi}
 \par Each natural number $n$ can be written in base $d\geq 2$, so as a finite word in the alphabet $\Sigma:=\{0,\ldots,d-1\}$. A subset of $\N$ is $d$-automatic if it is recognized by a finite automaton $\A$ working on $\Sigma$. 

\par Let $d=2$. The study of $2$-automatic subsets of $\N$ has been marked by a result of R. B\"uchi who showed that the sets definable in $(\N,+,V_2)$, where $V_2(n)$ is the highest power of $2$ dividing $n$, are exactly the $2$-automatic sets. (The original statement of R. B\"uchi is about $(\N,+,P_2)$, where $P_2$ is the set of powers of $2$, which is incorrect as observed, for instance, by Semenov. 
 In his review of B\"uchi's paper, McNaughton suggested to consider instead the structure $(\N,+,\epsilon_2)$, where $\epsilon_2(n,m)$ holds iff $n$ is a power of $2$ which occurs in the binary expansion of $m$. This last structure is interdefinable with $(\N,+,V_2)$ (a modified proof of B\"uchi's result with $V_2$ can be found in \cite{Br}).
  
\par Note that  $(\Z,+,<,V_2)$ has IP (one can define any finite subset of $P_2$ using $\epsilon_2$).
But the expansion $(\Z,+,<,P_2)$ is NIP \cite[Corollary 2.33]{LP}, \cite[Theorem 6.1]{H1}. This can be shown using the following model-theoretic description of definable sets or a criterium due to Chernikov and Simon on NIP pairs of structures.

\fct \cite{D} The theory of $(\Z,+,<,P_2)$ is model-complete in the language $\{+,-,<,\mod_n,\lambda_2\colon n\in \N\setminus\{0,1\}\}$, where $\lambda_2(x)$ is the largest power of $2$ smaller than $x$.
\efct
\par If one forgets about the order and considers the reduct: $Th(\Z,+,P_2)$, we have the following result. (We postpone till the next section, the definition of elementary $2$-sets and more generally of $F$-sets, but we recall here the notion of a stable subset.)

\dfn \cite[Definition 4.1]{H}, \cite{H1} A subset $A\subset \Z$ is stable if the formula $x+y\in A$ is stable in $Th(\Z,+,A)$.
\edfn
 \fct \label{stableZ} \cite{MS} The theory of $Th(\Z,+,P_2)$ is stable.
\par\noindent \cite{H} Let $A\subset \Z$. Then $A$ is $2$-automatic and stable in $(\Z,+)$ iff $A$ is a boolean combination of elementary $2$-sets and cosets of subgroups of $\Z$ iff $A$ is definable in $(\Z,+,P_2)$.
 \efct
 \section{Finitely generated abelian groups with a distinguished injective endomorphism}\label{main}
 \par In their work on Mordell-Lang conjecture for semiabelian varieties over finite fields, R. Moosa and T. Scanlon studied certain finitely generated abelian groups $\Gamma$ endowed with an endomorphism $F$ \cite[Section 2]{MS} and equipped them with a collection of sets called $F$-sets, for instance orbits by powers of $F$. They called the resulting structure on the abelian group $\Gamma$, an $F$-structure and showed it admits quantifier elimination and is stable (see \cite[Theorem A]{MS} for a precise statement). 
\par  Let $(\Gamma,+,0)$ be a finitely generated $\Z[F]$-module (namely a $\Z[X]$-module with the action of $X$ given by an endomorphism $F$ of $(\Gamma,+,0)$). 
We further assume that $F$ is injective.
 \par Note that the structure $(\Gamma,+,F)$, adding a function symbol for $F$, or $(\Gamma,+,0)$ endowed with the usual $\Z[F]$-module language, namely adding scalar multiplication for each element of $\Z[F]$ are bi-interpretable. In particular the theory of $(\Gamma,+,F)$ is stable.  
\par Let $\Sigma$ be a finite symmetric generating subset of $\Gamma$ as a $\Z[F]$-module, containing $0$. These will be the assumptions on $(\Gamma,+,F)$ from now on.
\medskip
  \par We say that the finite word $\sigma:=a_0\cdots a_n$ represents $g\in \Gamma$ if $g=a_0+F(a_1)+\ldots+F^n(a_n)$. We use the notation $g=[\sigma]_F$. 
  \par Here we adopt the definition of $F$-spanning subset given in \cite{H} which is slightly more general that the one given in  \cite[Definition 5.1]{BM}, or in \cite[Definition 2.4]{H2}.
  \dfn \cite[Definition 2.12]{H} \label{span} Let $\Sigma$ be as above.
  Then $\Sigma$ is a $F$-spanning set for $\Gamma$ if
 \begin{enumerate}[label = {(C\arabic*)}]
 \item any element of $\Gamma$ can be represented by an element of $\Sigma^*$,
 \item if $a_1, a_2, a_3\in \Sigma$, then $a_1+a_2+a_3\in \Sigma+F(\Sigma)$,
 \item if $a_1, a_2\in \Sigma$ are such that for some $b\in \Gamma$, $a_1+a_2=F(b)$, then $b\in \Sigma$.
 \end{enumerate}
 \edfn
\dfn  \cite[Definition 6.4]{BM}, \cite[Definition 2.5]{H2}  A subset $A\subset \Gamma$ is $F$-automatic if there is an $m>0$ and an $F^m$-spanning set $\Sigma$ such that 
$\{\sigma\in \Sigma^*\colon [\sigma]_F\in A\}$ is a regular language. (One easily extends the definition to subsets of $\Gamma^n$, $n>1$.)
\edfn
\par Then one shows that if a subset of $\Gamma$ is $F$-automatic, then for any $n>0$ and any $F^n$-spanning subset $\Sigma$, $\{\sigma\in \Sigma^*\colon [\sigma]_{F^n}\in A\}$ is a regular language \cite[Proposition 2.6]{H2}.
\medskip
\par When $\Gamma$ is a finitely generated group, C. Hawthorne has shown the following characterization of the existence of a spanning subset. Let $F$ an injective endomorphism of $\Gamma$. Then for some $n>0$ there is an $F^n$-spanning set for $\Gamma$ if and only if the eigenvalues of $F\otimes (1\restriction \C)$ in the $\C$-vector space $\Gamma\otimes \C$ are of modulus $>1$ \cite[Theorem 3.12]{H2}. 
\par In the general case, C. Hawthorne gave another characterization in terms of existence of {\it well-behaved} length functions \cite[Theorem 2.43]{H} (that he will use to show that $F$-sets are $F$-automatic (see Fact \ref{rec})).
\rem Note that in the definition of automatic structure, one chooses a set of representatives for the elements of the domain of the structure (and similarly for the predicates). In the definition above for a subset $A$ of $\Gamma$ to be $F$-automatic, one needs to check that all the finite words which represents an element of $A$ form a regular language. However, we have the following result of Bell and Moosa, which says that if $L$ is a regular language on $\Sigma$, then the subset $[L]_F$ of $\Gamma$ consisting of $\{[\sigma]_F\colon \sigma\in L\}$ is $F$-automatic \cite[Proposition 6.8 (2)]{BM}. This is a consequence of the fact that the equivalence relation $E$ defined between two words $\sigma, \tau$ (of the same length) by $E(\sigma,\tau)$ if and only if $[\sigma]_{F}=[\tau]_{F}$, is regular.
\par We have seen we may represent an element $g\in \Gamma$ by words of different length, but we use the extra symbol $\sharp$ and the operation of convolution.
\erem
 \dfn \label{Fsets}\cite[Definition 2.3]{MS}  Let $a\in \Gamma$, then set $K(a,F):=\{a+F(a)+\cdots+F^n(a)\colon n\in \N\}$. 

 \par An {\it elementary $F$-set} of $\Gamma$ is a subset of $\Gamma$ of the form $a_0+K(a_1,F^{n_1})+\ldots+K(a_m,F^{n_m})$, for some $a_0,\ldots,a_m\in \Gamma$ and $n_1,\ldots,n_m\in \N$. 
For instance, the $F$-orbit of $a$: $\{F^n(a): n\in \N\}=a+K(Fa-a,F)$ is an elementary $F$-set.
\par An {\it $F$-set} of $\Gamma$ is a finite union of sets which can be written as a sum of an elementary $F$-set of $\Gamma$ and an $F$-invariant subgroup of $\Gamma$ (namely $\Z[F]$-submodules).
\edfn
\par Let $\cF(\Gamma)$ (respectively $\cF(\Gamma^n)$) be the collection of $F$-sets of $\Gamma$ (respectively $\Gamma^n$). Consider the structure whose domain is $\Gamma$ and a set of predicates interpreted by the elements of $\bigcup_{n>0} \cF(\Gamma^n)$.
Denote the resulting structure by $(\Gamma, \cF)$. (Note that it expands $(\Gamma,+,F)$.) Assume $\Gamma$ is a finitely generated abelian group and $\bigcap_{i\in \N} (F^i)=\{0\}$, where $(F^i)$ denotes the ideal generated by $F^i$ in $\Z[F]$. Then Moosa and Scanlon showed that the theory of $(\Gamma,\cF)$ is stable \cite[Theorem 6.11]{MS} (see also \cite[Remark 6.12]{MS} where a proof of superstability is sketched).

\par Now we can state the result analogous to Fact \ref{stableZ}, for $\Z[F]$-module $(\Gamma,+,F)$, where $\Gamma$ is a finitely generated, but first we need to recall the notion of {\it sparse}. As before, a subset $A$ of $\Gamma$ is stable if the formula $x+y\in A$ is a stable formula.
\dfn Let $L\subset \Sigma^*$, then $L$ is {\it sparse} if $L$ is regular and if the set of words in $L$ of length smaller than or equal to $n$, $n\in \N$, is bounded by a polynomial function of $n$.  (For instance the set of binary expansions of powers of $2$ is sparse).
\par A subset $A\subset \Gamma$ is $F$-sparse if $A=[L]_{F^r}$ for some sparse $L\subset \Sigma^*$ with $\Sigma$ a $F^r$-spanning set, $r>0$.
\edfn
 \par  We denote the length of $\sigma\in \Sigma^*$ by $\ell(\sigma)$. Let $\ell_{\Sigma}:\Gamma\to \N$ be the map sending $g\in \Gamma$ to the length of the shortest word $\sigma=a_0\ldots a_n\in \Sigma^*$ such that $[\sigma]_F=g$.  
 We will say that a word $\sigma$ has support $1$ if there is only one occurrence in $\sigma$ of elements of $\Sigma\setminus\{0\}$.
 \fct \cite{MS}, \cite[Theorem 6.3]{H2} Assume that $\Gamma$ is a finitely generated group and that it has a $F^m$-spanning subset for some $m>0$. Let $A\subset \Gamma$ be $F$-sparse and stable in $\Gamma$. Then $A$ is a boolean combination of elementary $F$-sets.  
 \efct

\par Now starting with the following result of Bell and Moosa (\cite[Lemma 6.7, Theorem 6.9]{BM}), revisited by C. Hawthorne, we aim to define an expansion of $(\Gamma,+,F)$ where the definable subsets are exactly the $F$-automatic ones.
\fct \label{rec} \cite[Proposition 2.32, Theorem 2.54]{H} Let $\Gamma$ be a finitely generated $\Z[F]$-module and suppose it has a $F^m$-spanning set for some $m>0$. Then $(\Gamma,+,F)$ is FA-presentable.  Moreover every $F$-set of $\Gamma^m$, $m>0$, is $F$-automatic. 
\efct
\medskip
\par First let us make a few observation on $\Sigma$.
Note that if $a\in \Sigma\cap F(\Gamma)$, then by condition (C3), $a=F(b)$ for some $b\in \Sigma$. So w.l.o.g. we may assume that $\Sigma\subset \Gamma\setminus F(\Gamma)\cup\{0\}$ {\bf and it will be our assumption from now on}.
Then we pick a subset $\Sigma_0$ of $\Sigma$ containing $0$ and a unique coset representative of $\Sigma$ modulo $F(\Gamma)$. Then for $a_0\in \Sigma$, there is $a\in \Sigma_0$ such that $a_0-a\in F(\Gamma)$. So by condition $(C3)$, there is $b\in \Sigma$ such that $a_0+(-a)=F(b)$, namely $a_0=a+F(b)$. More generally, using this subset $\Sigma_0$, we want to represent any element of $\Gamma$ as follows.
\lem  \label{unicity} Let $g=a_0+\ldots+F^n(a_n)\in \Gamma\setminus\{0\}$ with $a_i\in \Sigma$, $1\leq i\leq n$, then there is another representation of $g$ of the form $F^m(b_m)+\ldots+F^{n}(b_n)+F^{n+1}(b)+F^{n+2}(b')$ with $b_m\in \Sigma_0\setminus\{0\}$, $m\leq i\leq n$, $b_i, b, b'\in \Sigma$ and where $m$ is the smallest indice where $a_i\neq 0$.
\elem
 \pr Let $g=a_0+F(a_1)+\ldots+F^n(a_n)\neq 0$. We will show the following statement: $g=F^m(b_m)+\ldots+F^{n}(b_n)+F^{n+1}(b)+F^{n+2}(b')$ with $b_m\in \Sigma_0\setminus\{0\}$, $m\leq i\leq n$, $b_i\in \Sigma_0, b, b'\in \Sigma$, where $m$ is the smallest indice where $a_i\neq 0$.
 \par By the choice of $\Sigma_0$, 
  there is $b_{0}\in \Sigma_0, a_1'\in \Sigma$ such that $a_0=b_0+F(a_1')$.
\par  Then we rewrite $g$ as $b_0+F(a_1+a_1')+F^2(a_2)+\ldots+F^n(a_n)$. By condition (C2), there are $u_1, u_2\in \Sigma$ such that $a_1+a_{1}'=u_{1}+F(u_{2})$. Now let $b_{1}\in \Sigma_{0}$ be such that $u_{1}=b_1+F(a_2')$ with $a_2'\in \Sigma$. 
\par Then we rewrite $g$ as $b_0+F(b_1)+F^2(a_2+u_{2}+a_{2}')+\ldots$. 
 Then write $a_2+u_2+a_2'=u_3+F(u_4)$, $u_3, u_4\in \Sigma$ and proceed as before, namely replace $u_3$ by $b_2+F(a_3')$ with $b_{2}\in \Sigma_{0}$. 
 So at stage $m\leq n$, we obtain $g=b_0+F(b_1)+\ldots+F^{m-1}(b_{m-1})+F^{m}(a_m+u+u')+ \ldots$, with $u, u'\in \Sigma$. Apply condition (C2), to write $a_m+u+u'$ as $v+F(w)$ for some $u, w\in \Sigma$ and $v$ as $b_{m}+F(v')$, $b_{m}\in \Sigma_{0}$. So we have replaced $F^{m}(a_m+u+u')$ by $F^m(b_{m})+F^{m+1}(w+v')$. Either $m=n$ and are left with a summand of the form $F^{n+1}(w+v')$ that we can write as $F^{n+1}(w')+F^{n+2}(w'')$, for some $w', w''\in \Sigma$, or $m<n$ and we have to deal an expression of the same form, namely $F^{m+1}(a_{m+1}+w+v')$.\qed
\medskip
\par From now on, we will also make the convention that when we write $g=[\sigma]_F$ with $\sigma\in \Sigma^*$, then the rightmost letters in $\sigma$ is not equal to $0$. 
For $\bc\in \Sigma\setminus\{0\}$, denote by $\Sigma^*\bc$ the subset of $\Sigma^*$ of all finite words ending with $\bc$. With this convention, we have that $\Gamma\setminus\{0\}=\bigcup_{\bc\in \Sigma\setminus\{0\}} [\Sigma^*\bc]_F$.
\medskip 
\nota Let $I_F:=\{F^n(a)\colon a\in \Sigma\setminus\{0\}, n\in \omega\}$. It is easily seen that this subset $I_{F}$ (of $\Gamma$) is automatic. Indeed in the $F$-representation of the elements of $\Gamma$, $I_{F}$ is the set of elements which can be represented by a (finite) word of support $1$. 
\par For each $a\in \Sigma\setminus\{0\}$, we will also need to single out the orbits $a$ by $F$: $I_{F,a}:=\{F^n(a)\colon n\in \omega\}$. This subset $I_{F,a}$ (of $\Gamma$) is again automatic and $I_{F}=\bigcup_{a\in \Sigma\setminus\{0\}} I_{F,a}$.
\enota
\par Recall that we assumed that $F$ is injective and that $\Sigma\setminus\{0\}\subset \Gamma\setminus F(\Gamma)$), so one cannot have $F^n(a)=F^{m}(b)$ with $a, b \in \Sigma\setminus\{0\}$ and $m\neq n$. This allows us to define the following preorder $\preceq$ on $I_F$.
\dfn We set $F^i(a)\preceq F^j(b)$ if $i\leq j$ and $F^i(a)\prec F^j(b)$ if $i<j$, with $a, b\in  \Sigma\setminus\{0\}$, (equivalently $F^i(a)\prec F^j(b)$ if ($F^i(a)\preceq F^j(b)$ and $\neg (F^j(b)\preceq F^i(a))$)). 
Furthermore we set $I_F\prec 0$.
\edfn
\par We will use the notation $u\sim v$, for $u, v\in I_{F}$ to mean that $u\preceq v$ and $v\preceq u$.
\par For convenience, we also define the partial function $F^{-1}$ on $I_F\setminus \Sigma$ as follows: let $u\in I_F\setminus \Sigma$, then $F^{-1}(u)=v \leftrightarrow u=F(v)$. 
\dfn For $g\in \Gamma\setminus\{0\}$, define a (unary) function $V_F:\Gamma\setminus\{0\}\to I_F$ as follows. 
We use Lemma \ref{unicity} to represent $g$ as $\sum_{i=m}^n F^i(b_i)$, with $b_m\neq 0$, $b_m\in \Sigma_0$, $b_i\in \Sigma$, $m\leq i\leq n$. Then define  $V_{F}(g):=F^m(b_{m})$. We extend $V_{F}$ on the whole of $\Gamma$ by $V_{F}(0)=0$.
\edfn
Note that $V_F$ is well-defined by choice of $\Sigma_0$ and that the graph of $V_F$ is automatic.
\par Let $g, g'\in \Gamma\setminus\{0\}$, $u\in I_{F}$. Then, using condition (C2), we have that if $u\prec V_{F}(g)$ and $u\prec V_{F}(g')$, then $u\prec V_{F}(g-g')$.
\medskip

 Also, one cannot have $F^n(a)=F^{m}(b)+F^k(c)$ with $a, b, c \in \Sigma\setminus\{0\}$ and $m<k$, and $m\neq n$, using that $\Sigma\cap F(\Gamma)=\{0\}.$
 \dfn We define the binary relation $R$ on $(\Gamma\setminus\{0\})\times (\Gamma\setminus\{0\})$ as follows $R(g,u)$ if $u\in I_F$ with $u=F^n(a)$, $a\in \Sigma\setminus\{0\}$ and $g=\sum_{i=m}^n F^i(a_i)$, with $a_{i}\in \Sigma$, $m\leq i\leq n$, $a_{m}\neq 0$ and $a_{n}=a$, namely $a$ occurs as the rightmost letter in $\sigma\in \Sigma^*$ with $g=[\sigma]_F$. (The relation $R(g,u)$ holds if there is $\sigma\in \Sigma^*$ such that $[\sigma]_F=g$ and the rightmost letter of $\sigma$ equals $u$. We can extend it to $\Gamma\times\Gamma$ by adding $R(0,0)$.)
\par  Then we define the relation $\epsilon_F(u,g)$, expressing that $u$ occurs in some representation of $g$, by
\[
u\in I_F\,\wedge\,(g=u\vee\;(\exists h\;(R(h,u)\wedge u\prec V_F(g-h))).
\]
\edfn

\par Note that both relations $R$ and $\epsilon_{F}$ are automatic. 

\lem  Let $g\in \Gamma$ and $u\in I_F$ and suppose that $\epsilon_F(u,g)$ holds. Then there is a unique element $h\in \Gamma$ such that
\[R(h,u)\wedge u\prec V_F(g-h).
\]
\elem

\pr Let $g, u$ as above with $u=F^m(a)$, $m\geq 0$, $a\in \Sigma$ and suppose that we have $h, h'$ such that $R(h,u)\wedge u\prec V_F(g-h)$ and $R(h',u)\wedge u\prec V_F(g-h')$. This implies that $u\prec V_F(h-h')$.
\par By definition, both $h, h'$ have a representation of the form $h=a_0+F(a_1)+\ldots+F^m(a)$, $h'=b_0+F(b_1)+\ldots+F^m(a)$, with $a_i, b_i\in \Sigma$, $0\leq i\leq m$. If $m=0$, then $h=h'$, so let us assume that $m\geq 1$.
Let $i$ be the smallest index such that $a_i-b_i\neq 0$. Write $a_i+(-b_i)=d+F(d')$, with $d, d'\in \Sigma$, applying condition (C2). If $d\neq 0$, we get a contradiction (with $u\prec V_F(h-h')$). So assume that $a_i+(-b_i)=F(d')$ with $d'\neq 0$ and replace $a_i$ by $b_i+F(d')$. If $i+1=m$, then $V_F(h-h')\sim u$, a contradiction. If $i+1<m$, we compare $a_{i+1}+d'$ and $b_{i+1}$. By condition $(C2)$, we have $a_{i+1}+d'+(-b_{i+1})=v+F(w)$, $v, w\in \Sigma$. If $v\neq 0$, again we get a contradiction. If $v=0$, then we replace $a_{i+1}+d'$ by $b_{i+1}+F(w)$. If $w=0$, then we have obtained another expression for $h, h'$ which makes them equal up to index $i+1$. If there is another index $m>j>i+1$ where they differ, we proceed as before. If $w\neq 0$ and $i+2<m$, then we compare $a_{i+2}+w$ and $b_{i+2}$, until we reach a contradiction or the case where we found a witness for $h=h'$. We eventually reach the case where $a_{i+k}=a$ and either we get that $h=h'$ or that $V_F(h-h')\sim u$, a contradiction. \qed
\rem Now suppose that $g\in \Gamma\setminus\{0\}$, $u\in I_F$ but $\neg \epsilon_F(u,g)$, then either there is a largest (in the preorder $\preceq$) $\tilde u\preceq u\in I_F$ such  that $\epsilon_F(\tilde u,g)$, then by the preceding lemma there is a unique element $h$ such that $R(h,\tilde u)\wedge \tilde u\prec V_F(g-h)$, otherwise set $h=0$. Note that in both cases, we still have that $u\prec V_F(g-h)$ and that we have a unique such $h$ with $\forall \tilde u\in I_F\;(R(h,\tilde u)\rightarrow \tilde u \preceq u)$.

\erem
\nota \label{interval} This allows us to define an auxiliary function $f:\Gamma\setminus\{0\}\times I_{F} \to \Gamma$ 
by $f(g,u)=h$ if ($V_F(g)\preceq u\;\wedge\;\exists u'\in I_{F}\;(u'\preceq u\,\wedge \,R(h,u')\,\wedge  u\prec V_F(g-h))$ or $u\prec V_F(g)\;\wedge h=0$).
\par\noindent  Let $u_1, u_2\in I_F$ with $u_1\prec u_2$. Then denote by 
 \[
 g\restriction [u_1\;u_2]:=f(g,u_2)-f(g,F^{-1}(u_1))
 \;{\rm  and}\; g\restriction [u_2\;u_1[:=f(g,u_2)-f(g,u_1).
 \]
\enota

\bigskip
\par Let $\L:=\{+, F, V_F, R,I_{F}, \preceq, I_{F,a}; a\in \Sigma\setminus\{0\}\}$ and 
\par\noindent 
$\Gamma_V:=(\Gamma,+,F,V_F, R, I_{F}, \preceq, I_{F,a}; a\in \Sigma\setminus\{0\}).$
\prop \label{F-def} Let $(\Gamma,+,F)$ be a finitely generated $\Z[F]$-module, where $F$ is injective, with a $F$-spanning subset. Then $\Gamma_{V}$ is FA-presentable. In particular this structure is decidable. The structure $(\Gamma,\cF)$ is definable in $\Gamma_V$. 
\eprop
\pr The fact that $\Gamma_V$ is FA-presentable follows from \cite[Proposition 2.32]{H} where it is shown that the diagonal, the graph of $+$ and the graph of $F$ are $F$-automatic. In the discussion above we also noted that all the predicates $V_F, R, I_{F}, \preceq, I_{F,a}; a\in \Sigma\setminus\{0\}$ are automatic. Now let us show that any elementary F-set is definable in $\Gamma_V$. 
Note that it will give another proof that $F$-sets are automatic \cite[Theorem 2.54]{H}.
\par First let us define the elementary F-set $K(a,F)$, $a\in \Sigma\setminus\{0\}$, as follows:
 $g\in K(a; F)$ if 
\[
\exists u\in I_{F,a}\;(R(g,u)\;\wedge\;
\forall c\in I_{F}(V_F(g)\preceq c\prec u\;\rightarrow g\restriction [c\;F(c)[\in I_{F,a})).
\]
\par Now let $a\in \Gamma\setminus\{0\}$, $a=\sum_{i=0}^m F^{n_i}(a_i)$, $a_i\in \Sigma\setminus\{0\}$, $0\leq n_{0}<\ldots<n_{m}$ and $k=n_{m}-n_{0}+1$. First define the orbit of $a$, namely $Orb(a):=\{F^m(a): m\in \N\}$ as follows: $\exists g_0\;\ldots\exists g_m\;$
\[
\bigwedge_{i=0}^m g_i\in F^{n_i}(I_{F,a_i})\;\wedge\;g=\sum_{i=0}^m g_i\,\wedge\;\bigwedge_{i=0}^{m-1} F^{n_{i+1}-n_i}(g_{i})\sim g_{i+1}
\]
Define $g\in K(a;F^k)$ as follows:
\[
\exists h\;\big(\exists u\in I_F\;(V_F(h)=V_F(g)\;\wedge\;R(g,u)\wedge R(h,F(u))\wedge\;
\]
\[
(\forall u_1\forall u_2 \big(u_1, u_{2}\in I_{F} \wedge u_{1}\prec u_{2} \wedge  \bigwedge_{i=1}^2\epsilon_{F}(h,u_{i})\wedge  \neg(\exists u_3\in I_F (\epsilon_{F}(h,u_{3})\wedge  u_1\prec u_3\prec u_2)))\big)
\]
\[
\rightarrow\;(g\restriction [u_1\;u_2[ \in Orb(F^k(a))\big)\big).
\] 
Finally, define $g\in K(a;F)$ as follows:
\[
\exists g_0\;\ldots\exists g_{k-1}\;
\bigwedge_{i=0}^{k-1} g_i\in K(F^i(a);F^k)\wedge g=\sum_{i=0}^{k-1} g_{i}.
\]
\par One can easily modify the above formulas to define any elementary F-set, as well as any elementary F-set of some power of $\Gamma$.
 \qed
\medskip
\par We have that as an $\Z[F]$-module, the structure $(\Gamma,+,F)$ is stable and as we have seen earlier on certain conditions on $(\Gamma,+,F)$, $(\Gamma,\cF)$ is stable (even superstable) \cite[Theorem 6.11]{MS}. 
\par C. Hawthorne considered the question of which subsets $A$ of $\Gamma$, to add to get a NIP-expansion $(\Gamma,+,A)$ \cite[section 7]{H2}.
\par In our framework, it would be interesting to find an intermediate structure between $(\Gamma,+,F)$ and $\Gamma_V$ which is NIP, more specifically:
\par {\bf Question} Is there a reduct of $\Gamma_V$ interdefinable  with $(\Gamma,+, F,\cF(\Gamma))$, as in the special case of $\Gamma=\Z$ and $F$ the multiplication by a prime number? 

\medskip
 \par Now we want to identify all automatic sets in some $\Gamma^n$, $n>0$, as the class of all $\L$-definable subsets in $\Gamma_{V}$.
The proof presented here follows the same strategy as in \cite{MP}.
We use Kleene's theorem on the correspondence between regular languages and languages accepted by a finite automaton (see section \ref{aut}). 
\par A well-known observation is that if a language $L$ is accepted by a finite automaton, then the language consisting of all words in $\Sigma^*$ of the form $\underbrace{0\ldots0}_n\omega$ with $\omega\in L$ is again recognizable. We will denoted such language by $0^nL$, $n>0$ and 
 when $n=0$, $0^{0}=\lambda$.

\lem \label{regular} Let $L$ be a regular language on $\Sigma$ and let $\bc\in \Sigma\setminus\{0\}$. Then there are $\L$-formulas such that  $\rho_{L}(g,x)$ such that
$\Gamma_V\models \rho_{L}([w^{\frown}\bc]_{F},F^n(\bc))$ iff $w\in 0^nL$. 
\elem
We work by induction on the complexity of regular languages. 
Let $L$, $L_1$, $L_2$ be languages of complexity $\leq c+1$. 
Since any element of $\Gamma\setminus\{0\}$ is represented by a word ending with an element of $\Sigma\setminus\{0\}$, we consider elements of $\Sigma^*$ concatenated with a letter in $\bc\in \Sigma\setminus\{0\}$. Moreover for each such letter $\bc$, there is some $b\in \Sigma$ such that $\bc=\bc_0+F(b)$ with $\bc_0\in \Sigma_0$.
\par By induction we assume that for languages $L$ of complexity $\leq c$, there exist formulas $\rho_{L}(g,x)$ such that
$\rho_{L}([w^{\frown}\bc]_{F},F^n(\bc))$ holds iff $w\in 0^nL$. 
\par $\bullet$ The case of complexity $0$ is the case of languages $L$  of the form $L=\{a\}$ with $a\in \Sigma$. 
\par \noindent If $L=\{0\}$, then set $\rho_{L}(g, x):=(g=F(x))$;\\
 if $L=\{a\}$, $a\neq 0$, then set $\rho_{L}(g, x):=(g=u+F(x)\wedge u\in I_{F,a}\wedge F(u)\sim F(x)).$ So $\rho_{L}([w^{\frown} \bc]_{F},F^n(\bc))$ holds iff 
$[w^{\frown} \bc]_F=F^n(a)+F^{n+1}(\bc)$ iff $w=0^na$.
\par $\bullet$ Then by induction we assume that we have formulas $\rho_{L_i}(g,x)$ such that $\rho_{L_i}([w^{\frown} \bc]_F,F^n(\bc))$ holds iff $w\in 0^n(L_i)$, $0\leq i\leq 2$ and we show there exist formulas with the following property:
\begin{enumerate}
\item $\rho_{L_1\cup L_2}([w^{\frown} \bc]_{F},F^n(\bc))$ holds iff $w\in 0^n(L_1\cup L_2)$, 
\item $\rho_{L_1^{\frown} L_2}([w^{\frown} \bc]_{F},F^n(\bc))$ holds iff $w\in 0^n(L_1^{\frown} L_2)$, $w\bc=0^nv_{1}v_{2}\bc, v_{1}\in L_{1}, v_{2}\in L_{2}$,
\item $\rho_{L_0^*}([w^{\frown} \bc]_{F},F^n(\bc))$ holds iff $w\in 0^n L_0^*$, $w\bc=0^nv_{1}v_{2}\ldots v_{m}\bc, v_{i}\in L_{0}$.
\end{enumerate}
\par Case $(1)$ is immediate since $w\in 0^n(L_1\cup L_2)$ iff ($w\in 0^nL_1$ or $w\in 0^nL_2$).
\par In case $(2)$, we let $\rho_{L_1^{\frown} L_2}(g,x)$ be the following formula: 
\[\exists v, u\in I_{F,\bc}\exists v', v''\in I_{F,\bc}\;((x \preceq v \preceq u\,\wedge \,S(v)\sim v'\,\wedge v''\sim S(u))\wedge  
\]
\[
x\preceq V_{F}(g)\,\wedge \,R(g,u)\wedge \,\rho_{L_1}(f(g,v)+v',x)\wedge \, \rho_{L_2}(g-f(g,v)+v'',v').
\]
Then we check that if $\rho_{L_1^{\frown} L_2}([w^{\frown} \bc]_F,F^n(\bc))$ holds iff $w\in 0^{n}(L_1^{\frown} L_2)$.
\par In case (3), 
we let $\rho_{L_0^*}(g,x)$ to be the formula:
\[
g=x\vee \exists h\;\exists u\in I_{F,\bc}\big(R(h,u)\, \wedge  R(g,u)\,\wedge  (V_F(h)\sim x\sim V_F(g))\,\wedge  
\]
\[
(\forall u_1\forall u_2 \big(u_1, u_{2}\in I_{F,\bc} \wedge u_{1}\prec u_{2} \wedge  \bigwedge_{i=1}^2\epsilon_{F, \bc}(h,u_{i})\wedge  \neg(\exists u_3\in I_F (\epsilon_{F}(h,u_{3})\wedge  u_1\prec u_3\prec u_2)))\big)
\]
\[
\rightarrow\; \rho_{L_{0}}(g\restriction [u_1\;u_2],u_1)\big)
\] 

Again, we check that $\rho_{L_0^*}([w^{\frown} \bc]_F,F^n(\bc))$ holds iff $w\in 0^{n}L_0^*$.\qed
\thm Let $(\Gamma,+,F)$ be a finitely generated $\Z[F]$-module, where $F$ is injective, with a $F$-spanning subset. Then $\Gamma_V$ is FA-presentable and definable subsets in $\Gamma_{V}$  
coincide with automatic sets.
\ethm
\pr The first part of the statement is Proposition \ref{F-def}. The second part follows from Lemma. \ref{regular}. 
We will simply do the proof for subsets of $\Gamma$. For subsets of $\Gamma^n$, one uses the operation of convolution (and works with the alphabet $(\Sigma\cup\{\sharp\})^n$) (see section \ref{automatic}). 
\par Let $\Sigma$ be a $F$-spanning set for $\Gamma$ and let $A$ be an automatic subset of $\Gamma$. So there is a regular language $L\subset \Sigma^*$ such that $A=\{[\omega]_{F}\colon \omega\in L\}$.
Now the set of finite words ending with a non zero element of $\Sigma$ is a regular language, say $\tilde L$. So $L\cap (\tilde L\cup \{\lambda\})$ is again a regular language, that we rename $L$.
Let $g\in A$; then $g=[\omega]_F$, for some $\omega\in L$. First assume that $g\neq 0$, then by Lemma \ref{regular}, $\omega\in L$ iff $\rho_{L}([\omega^{\frown}\bc]_{F},\bc)$ holds.
\par Now $[w^{\frown}\bc]_{F}=g+F(u)$ with $u\in I_{F,c}\wedge  R(g,u)$. If $g=0$, then $[\lambda^{\frown}\bc]_{F}=\bc$ and $\lambda\in L$ iff $\rho_L(\bc,\bc)$.
\par So, putting both cases together, we have $g\in A$ iff 
$$ \Gamma\models(g\neq 0\wedge  (\exists u\;(u\in I_{F,c}\wedge  R(g,u)\wedge   \rho_{L}(g+F(u),\bc))\vee (g=0\wedge  \rho_{L}(\bc,\bc))).$$  \qed
\lem The structure $\Gamma_V$ has the independence property. When $\Gamma$ is a finitely generated abelian group and $\bigcap_{i\in \N} (F^i)=\{0\}$, $\Gamma_V$ is a proper expansion of $(\Gamma,\cF)$. 
\elem
\pr The structure $\Gamma_V$ has IP. As in the case of $(\Z,+,V_{2})$, one codes finite subsets of $I_{F}$ as follows.
Let $u\in I_{F}$ and $g\in \Gamma\setminus\{0\}$. Then $u\in g$ if $\epsilon_{F}(u,g)$ holds. 
\par When $\Gamma$ is a finitely generated $\Z$-module and $\bigcap_{i\in \N} (F^i)=\{0\}$, $(\Gamma,\cF)$ is stable \cite[Theorem 6.11]{MS}.
\qed
\subsection*{Acknowledgements} I would like to thank P. D'Aquino who gave me the opportunity to give a talk on that subject in the model theory  seminar in may 2022 in the University ''degli studi della Campania Luigi Vanvitelli''. I also would like to thank R. Moosa and C. Hawthorne for helpful communications.

\end{document}